# REGIONAL EXPONENTIAL GENERAL OBSERVABILITY IN DIFFUSION NEUMANN PROBLEM

## Zainab A. Jaafar[1],   Raheam A. Al-Saphory[2*]

[1]Department of Mathematics, College of Education / Tuz Khurmatu, Tikrit University, Tikrit, IRAQ

[1]Zainabali69@hotmail.com

[2]Department of Mathematics, College of Education for Pure Sciences, Tikrit University, Tikrit, IRAQ.

[2]Email:  saphory@hotmail.com

*To whom correspondence should be addressed

**KYWORDS:** Strategic sensors, $\omega$-general strategic sensors, $\omega_{E_G}$-observability, Neumann boundary conditions.

**ABSTRACT**: The aim of this paper is to extend the concept of regional exponential general observability to the case of Neumann boundary conditions problem in diffusion system. More precisely, for linear distributed parameter diffusion systems, we show that the sensors characterizations allow to reconstruct the regional exponential state in a sub-region of the considered systems domain. Moreover, various interesting results associated with the choice of sensors are given and illustrated in specific situations in order regional exponential general observability notion to be achieved. Finally, we also show that, there exists a dynamical system for diffusion system is not exponential general observable in the usual sense, but it may be regional exponential general observable.

## 1. ITRODUCTION

The exponential observation problems in usual sense of distributed parameter systems has received much attention in ref.s [1-2] and references therein. Recently, this notion was extended  to the regional analysis by Al-Saphory, EI Jai [3-5] and Zerrik [6] this notion developed to many concepts in finite time interval or in infinite [7-8]. One of the interested concepts in regional asymptotic analysis, is the asymptotic observability notion which was introduced and developed by Al-Saphory and El Jai  [9-13].

These concepts  depend on the observation of the initial state on a given sub-region $\omega$ of the system domain $\Omega$. The main reason behind the study of this work there exist some problems in the real world cannot observe the system state in the whole domain, but in a part of this domain [9,14]. The scenario described by (Figure 1) below, one is interested in reconstruction the state in the green zone $\Omega_C$ rather than in the entire space [14].

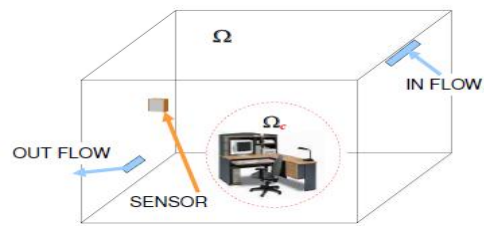

**Figure 1:** Observation region $\Omega_C$  with sensor

This problem is represented as regional exponential observation with interested  regional state in $\Omega_C$ of the domain $\Omega$. The purpose of this paper is to give an extension of the previous work  as in ref.s [15-16] to case of  Neumann boundary conditions by using strategic internal sensors. The paper is organized as follow. Section 2 gives an introduction of exponential regional  general observability in diffusion system problems. Section 3 deals with the main results. Finally, the last section illustrates applications of diffusion systems in many situations of sensors locations.

## 2. DIFFUSION SYSTEM OBSERVATION

### 2.1 Problem Statement

Let $\Omega$ be an open bounded set of  $R^n$, with smooth boundary $\partial\Omega$. Let $[\cdot, T]$, for $T > \cdot$ be a time measurement interval and $\omega$ be a nonempty given sub-region of $\Omega$. We denote $Q = \Omega \times ]\cdot, \infty[$ and  $\Sigma = \partial\Omega \times ]\cdot, \infty[$. Assume that $Z$, $U$, and $\mathcal{O}$ be separable Hilbert space where $Z = H^{\backslash}(\Omega)$ is the state space $U$ is the control space and $\mathcal{O}$ is the observation space. Usually, we consider $U = L^{\mathsf{Y}}(0, \infty, R^p)$ and $\mathcal{O} = L^{\mathsf{Y}}(0, \infty, R^q)$, where $p$ and $q$ hold for the number of  actuators and sensors [1-2].We consider the problem described by the following parabolic diffusion systems

$$\begin{cases} \frac{\partial z}{\partial t}(\mu, t) = \Delta z(\mu, t) + Bu(t) & Q \\ \frac{\partial z}{\partial v}(\eta, t) = \cdot & \Sigma \\ z(\mu, \cdot) = z_{\cdot}(\mu) & \Omega \end{cases} \qquad (1)$$

Augmented with the output function

$$y(., t) = Cz(., t) \qquad (2)$$

where $z_{\cdot}(\mu)$ supposed to be in $H^{\backslash}(\Omega)$ and unknown. The system (1) is defined with a Neumann boundary conditions, $\frac{\partial z}{\partial v}$ holds for outward normal derivative [17] and the operator generates $\Delta$ a strongly continuous semi-group $(S_{\Delta}(t))_{t \geq}$ on $Z$ [2] and is self-adjoint with compact resolvent. The operators $B \in \mathcal{L}(R^p, Z)$ and $C \in \mathcal{L}(Z, R^q)$ depend on the structures of actuators and sensors [1, 18], see (Figure 2).

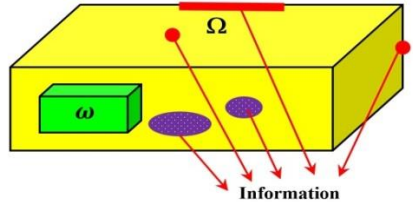

**Figure 2:** The domain of $\Omega$, the sub-region on $\omega$, the sensors locations

Under the given assumption, the system (1) has a unique solution [13],

$$z(\mu, t) = S_{\Delta}(t)z_{\cdot}(\mu) + \int_{\cdot}^{t} S_{\Delta}(t - \tau)Bu(\tau)\,d\tau \qquad (3)$$

The problem is that, how to reconstruct exponentially the current state in a given sub-region $\omega$, using internal zone or pointwise strategic.

## 2.2 $\omega$-strategic sensor

This sub-section is to extend sensors characterization notion in ref.[3,15-16], in order that, the diffusion system (1) is regionally exponentially general observable in Hilbert space $H^{\backslash}(\omega)$.

● Sensors are any couple $(D_i, f_i)_{\backslash \leq i \leq q}$ where $D_i$ denote closed subsets of $\bar{\Omega}$, which is spatial supports of sensors and $f_i \in H^{\backslash}(D_i)$ define the spatial distributions of measurements on $(D_i)$.

●The measurements may be given by $q$ sensors $(D_i, f_i)_{\backslash \leq i \leq q}$ and then, the output functions are given by the form:

$$y(., t) = \left[ y_{\backslash}(., t), y_{\mathsf{Y}}(., t), \dots, y_q(., t) \right]^{t_{\mathsf{Y}}}$$

where $\left[ y_{\backslash}(., t), y_{\mathsf{Y}}(., t), \dots, y_q(., t) \right]^{t_{\mathsf{Y}}}$ is the transpose of $y(., t)$ and then

$$y(., t) = Cz(., t) = \int_{D_i} z(\mu, t) f_i(\mu)\,d\mu = \langle z(., t), f_i(.)) \rangle_{H^{\backslash}(\Omega)} \qquad (4)$$

when $D_i \subset \Omega$.Thus, in the case of pointwise or filament we have

$$y(., t) = Cz(., t) = \int_{\Omega} z(\mu, t)\,\delta_{b_i}(\mu - b_i)\,d\mu = \langle z(., t), \delta_{b_i}(.) \rangle_{H^{\backslash}(\Omega)} \quad (5)$$

● In the case of internal pointwise sensors and filament, the operator $C$ is unbounded and some precautions must be taken in [1].

● For $z(\mu, t) = S_A(t)z,(\mu)$, defines the operator

$$K: Z \to L^\gamma(0, \infty, R^q)$$
$$z \to Kz = CS_A(t)z \qquad (6)$$

then $y(., t) = K(t)z,.$. Thus, in the case of internal zone sensors, linear and bounded [7]. The adjoint operator $K^*$ of $K$ is defined by

$$K^* y(., t) = \int_0^t S^*_A(s) \, C^* y(\xi, s) ds \qquad (7)$$

● For the region $\omega$ of the domain $\Omega$, the restriction operator $\chi_\omega$ is defined by

$$\chi_\omega : H^\gamma(\Omega) \longrightarrow H^\gamma(\omega)$$
$$z \longrightarrow \chi_\omega z = z \mid_\omega \qquad (8)$$

where $z \mid_\omega$ is the restriction of the state $z$ to $\omega$.

**Definition 2.1:** An autonomous diffusion system associated to (1)-(2) is exactly (respectively weakly) $\omega$-observable if:

$$Im \, \chi_\omega K^* = H^\gamma(\omega) \qquad \text{(respectively } \overline{Im \, \chi_\omega K^*} = H^\gamma(\omega).\text{).}$$

**Definition 2.2:** A sequence of sensors $(D_i, f_i)_{\gamma \leq i \leq q}$ are $\omega$-strategic if the diffusion systems (1)-(2) are weakly $\omega$-observable .

The concept of $\omega$-strategic sensor has been extended to the regional boundary case as in [19]. Assume that the set $(\varphi_{nj})$ of eigenfunctions of $H^\gamma(\Omega)$ orthonormal in $H^\gamma(\omega)$ associated with eigenvalues $\lambda_n$ of multiplicity $r_n$ and suppose that the system (1) has $J$ unstable modes. Then we have the following result.

**Proposition 2.3:** For the diffusion system (1), let $\sup r_n = r$, then suite of sensors $(D_i, f_i)_{\gamma \leq i \leq q}$ are $\omega$-strategic if and only if :

1. $q \geq r$

2. rank $G_n = r_n$, $\forall n$, $n = \gamma, ..., J$ with

$$G_n = (G_n)_{ij} = \begin{cases} (\frac{\partial \varphi_{nj}}{\partial v}, f_i(.))_{L^\gamma(D_i)} & \text{in the zone case} \\ \frac{\partial \varphi_{nj}}{\partial v} (b_i) & \text{in the pointwise case} \end{cases} \qquad (9)$$

**Proof:** The proof of this proposition is similar to the rank condition in ref. [5], the main difference is that the rank condition

rank $G_n = r_n$, $\forall n$, $n = \gamma, ...,$

for the proposition (2.1) need only hold for rank $G_n = r_n$, $\forall n$, $n = \gamma, ..., J$.

That means, when the rank is held for every $n = \gamma, \Upsilon, ...,$ then, obviously held for $n = \gamma, ..., J$ where $J$ is finite number.
∎

# 3. DIFINITIONS AND CHARACTERIZATIONS

Regional exponential general observability concept in diffusion system needs to extend some notions of an exponential behaviour (stability- detectability- estimator) in Hilbert space $H^\gamma(\omega)$ as ref.s [12-13].

**Definition 3.1:**A Semi-group is regionally exponentially general stable by the operator $T_\omega\Delta$ in $H^`(\omega)$  or ( $\omega_{E_G}$-stable) if, for every initial state $z_.(.) \in H^`(\Omega)$, the solution of the autonomous system associated with (1) converges exponentially to zero when $t \to \infty$, where $T \in H^`(\Omega)$ an operator with $T_\omega = \chi_\omega T$, is defined by $y(\mu, t) = T_\omega z(\mu, t)$.

**Definition 3.2:**The diffusion system (1) is said to be exponential general stable on $H^`(\omega)$ or ($\omega_{E_G}$-stable) if the operator $T_\omega\Delta$ generates a semi-group which is ($\omega_{E_G}$-stable). It is easy to see that the diffusion system (1) is $\omega_{E_G}$-stable if and only if, for some positive constants $F_\omega$ and $\sigma_\omega$,

$$\left\|S_{T_\omega\Delta}(.)\right\|_{\mathcal{L}(Z,\ H^`(\omega))} \le F_\omega e^{-\sigma_\omega t} \ t \ge \cdot \tag{10}$$

If $(S_{T_\omega\Delta}(t))_{t\ge.}$ is $\omega_{E_G}$-stable, then, for all $z_.(.) \in H^`(\Omega)$, the solution of autonomous diffusion system associated with (1) satisfies

$$\|z(t)\|_{H^`(\omega)} =\| S_{T_\omega\Delta}(.)z_. \|_{H^`(\omega)} \le F_\omega e^{-\sigma_\omega t}\|z_. \|_{H^`(\omega)} \tag{11}$$

and then

$$\lim_{t\to\infty}\|z(t)\|_{H^`(\omega)} = \cdot \tag{12}$$

**Definition 3.3:** The diffusion systems (1)-(2) are said to be exponential general detectable in $H^`(\omega)$or ( $\omega_{E_G}$-detectable) if there exists an operator $\mathcal{H}_\omega : R^q \to H^`(\omega)$ such that $(T_\omega\Delta - \mathcal{H}_\omega C)$ generates a strongly continuous semi-group $(S_{T_\omega L_\omega}(t))_{t\ge.}$ which is $\omega_{E_G}$-stable.

**Definition 3.4:** For the diffusion systems (1)-(2) consider the dynamical system

$$\begin{cases} \frac{\partial w}{\partial t}(\mu, t) = L_\omega w(\mu, t) + G_\omega u(t) + \mathcal{H}_\omega y(\mu, t) & Q \\ \frac{\partial w}{\partial v}(\eta, t) = \cdot & \Sigma \\ w(\mu, \cdot) = w_.(\mu) & \Omega \end{cases} \tag{13}$$

where $L_\omega$ generates a strongly continuous semi-group $(S_{L_\omega}(t))_{t\ge.}$ which is stable on Hilbert space W, $G_\omega \in \mathcal{L}(R^p, W)$ and $\mathcal{H}_\omega \in \mathcal{L}(R^q, W)$.The diffusion system (13) defines an $\omega_{E_G}$-estimator  for  $\chi_\omega T z(\mu, t)$ if

(1) $\lim_{t\to\infty} \|\chi_\omega T z(\mu, t) - w(\mu, t)\|_{L^`(\omega)} = \cdot$

(2) $\chi_\omega T$ maps $D(\Delta)$ in $D(L_\omega)$ where $z(\mu, t)$ and $w(\mu, t)$ are  the

  solution of the diffusion systems (1) and (13).

**Definition 3.5:**The diffusion systems (1)-(2) are regional exponential general observable in $\omega$ ($\omega_{E_G}$- observable) if there exists a dynamical system (13) ($\omega_{E_G}$-estimator) for the original system.

We can extend definitions 2.1 and 2.2 to the case of general observability for regional general decomposed diffusion systems [16],where $Z = z_\backslash \oplus z_\Upsilon$ and $Z = z_\backslash + z_\Upsilon$ with $Z \in D(T_\omega\Delta), Z_\backslash \in D(T_\omega\Delta_\backslash), Z_\Upsilon \in D(T_\omega\Delta_\Upsilon)$ as in [16].

**Definition 3.6:**The diffusion systems (1)-(2) are exactly (respectively weakly) general observable or ($\omega_G$-observable ) if:

$$Im\ \chi_\omega K^* = Z_\backslash \quad \text{(respectively } \overline{Im\ \chi_\omega K^*} = Z_\backslash \text{).}$$

**Definition 3.7:**A suite of sensors $(D_i, f_i)_{\backslash \le i \le q}$ are $\omega$ general strategic or ($\omega_G$-strategic) if the diffusion systems (1)-(2) are weakly $\omega_G$-observable.

  Thus, in this section is related to, present the sufficient condition of regional exponential general observability concept in diffusion system.

**Theorem 3.8:** The diffusion systems (1)-(2) are $\omega_{\varepsilon_c}$-observable by the dynamical system (13), if the following conditions hold:

(1) There exist $M \in \mathcal{L}(R^q, H^{\backslash}(\omega))$ and $N \in \mathcal{L}(H^{\backslash}(\omega))$ such that

$$MC + NT_\omega = I \qquad (14)$$

(2) $\begin{cases} T_\omega \Delta - L_\omega T_\omega = \mathcal{H}_\omega C \\ and \quad G_\omega = T_\omega B \end{cases}$ (15)

**Proof:** For $x(\mu, t) = T_\omega z(\mu, t)$ and $w(\mu, t)$ solution of (13), denote

$e(\mu, t) = x(\mu, t) - w(\mu, t)$.

We have

$\begin{aligned}
\frac{\partial e}{\partial t}(\mu, t) &= \frac{\partial x}{\partial t}(\mu, t) - \frac{\partial w}{\partial t}(\mu, t) \\
&= T_\omega \Delta z(\mu, t) + T_\omega B u(t) \\
&\quad -L_\omega w(\mu, t) - G_\omega u(t) - \mathcal{H}_\omega y(\mu, t) \\
&= L_\omega e(\mu, t) - L_\omega x(\mu, t) + T_\omega \Delta z(\mu, t) \\
&\quad -\mathcal{H}_\omega y(\mu, t) + T_\omega B u(t) - G_\omega u(t) \\
&= L_\omega e(\mu, t) + [T_\omega \Delta z(\mu, t) - L_\omega T_\omega z(\mu, t) \\
&\quad -\mathcal{H}_\omega C z(\mu, t)] + [T_\omega B - G_\omega] u(t) \\
&= L_\omega e(\mu, t) + [T_\omega \Delta - L_\omega T_\omega - \mathcal{H}_\omega C] z(\mu, t) \\
&\quad +[T_\omega B - G_\omega] u(t) = L_\omega e(\mu, t)
\end{aligned}$

Consequently, $e(\mu, t) = S_{L_\omega T_\omega}(t) e(\cdot, t)$ where $e(\cdot, t) = T_\omega z.(\mu) - w.(\mu)$. Since the operator $L_\omega T_\omega$ generates a strongly continuous semi-group $(S_{L_\omega T_\omega}(t))_{t\geq}$. which is $\omega_{\varepsilon_c}$-stable, then by using equation (10) we obtain

$$\left\| e(\mu, t) \right\|_{H^{\backslash}(\omega)} \leq F_L \ e^{-\sigma_L t} \left\| e.(t) \right\|_{H^{\backslash}(\omega)},$$

and therefore $\lim_{t\to\infty} e(\mu, t) = \cdot$.

Consider now $\hat{z}(\mu, t) = My(\mu, t) + Nw(\mu, t)$, then we have

$\begin{aligned}
\tilde{z}(\mu, t) &= z(\mu, t) - \hat{z}(\mu, t) \\
&= z(\mu, t) - My(\mu, t) - Nw(\mu, t) \\
&= z(\mu, t) - MCz(\mu, t) - NT_\omega z(\mu, t) \\
&\quad +N[T_\omega z(\mu, t) - w(\mu, t)] \\
&= N[T_\omega z(\mu, t) - w(\mu, t)] \\
&= N[x(\mu, t) - w(\mu, t)] \\
&= Ne(\mu, t).
\end{aligned}$

Consequently we have $\lim_{t\to\infty} \left\| \hat{z}(\mu, t) - z(\mu, t) \right\|_{H^{\backslash}(\omega)} = \cdot$. Thus, we have $\lim_{t\to\infty} \left\| T_\omega z(\mu, t) - w(\mu, t) \right\|_{H^{\backslash}(\omega)} = \cdot$ and the dynamical system (13) is $\omega_{\varepsilon_c}$-estimator. Finally the diffusion systems (1)-(2) are $\omega_{\varepsilon_c}$-observable.

∎

**Remark 3.9:** From the previous theorem (3.1) we can deduce the following results:

**1.** Conditions (14) and (15) in theorem (2.10) guarantee that the diffusion systems (1)-( 2) are $\omega_{\varepsilon_c}$-observable.

**2.** If a diffusion system which is $\Omega_{\varepsilon_c}$-observable, then it is $\omega_{\varepsilon_c}$-observable.

**3.** If a diffusion system is $\omega_{E_c}$-observable in $\omega$, then it is $\omega'_{E_c}$-observable. In every subset $\omega'$ of $\omega$, but the converse is not true. This may be proven in the following example:

**Example 3.10:** Consider the diffusion system:

$$\begin{cases} \frac{\partial z}{\partial t}(\mu_\backslash,\mu_\tau,t) = \frac{\partial^\tau z}{\partial \mu_\backslash^\tau}(\mu_\backslash,\mu_\tau,t) + \frac{\partial^\tau z}{\partial \mu_\tau^\tau}(\mu_\backslash,\mu_\tau,t) & Q \\ \frac{\partial z}{\partial v}(\eta_\backslash,\eta_\tau,t) = \cdot & \Sigma \\ z(\mu_\backslash,\mu_\tau,\cdot) = z_\cdot(\mu_\backslash,\mu_\tau) & \Omega \end{cases} \qquad (16)$$

and suppose there is a single sensor $(b_i,\delta_{b_i})$ with $D_i = \{b_i\}$. The output function

$$y(t) = \int_\Omega z(\mu_\backslash,\mu_\tau,t)\,\delta(\mu_\backslash - b_\backslash,\mu_\tau - b_\tau)d\mu_\backslash\, d\mu_\tau \qquad (17)$$

Where $\Omega = ]\cdot,a_\backslash[\times[\cdot,a_\tau[$ and $b_i \in \Omega$ are the location of sensors $(b_i,\delta_{b_i})$. The operator $\Delta = (\frac{\partial^\tau}{\partial \mu_\backslash^\tau} + \frac{\partial^\tau}{\partial \mu_\tau^\tau})$ generates a strongly continuous semi-group $(S_\Delta(t))_{t\ge\cdot}$ on Hilbert space $H^\backslash(\Omega)$ [20]. Consider the dynamical system

$$\begin{cases} \frac{\partial w}{\partial t}(\mu_\backslash,\mu_\tau,t) = \frac{\partial^\tau w}{\partial \mu_\backslash^\tau}(\mu_\backslash,\mu_\tau,t) + \frac{\partial^\tau w}{\partial \mu_\tau^\tau}(\mu_\backslash,\mu_\tau,t) & Q \\ \qquad\qquad + \mathcal{H}C\big(w(\mu_\backslash,\mu_\tau,t) - z(\mu_\backslash,\mu_\tau,t)\big) \\ \frac{\partial w}{\partial v}(\eta_\backslash,\eta_\tau,t) = \cdot & \Sigma \\ w(\mu_\backslash,\mu_\tau,\cdot) = w_\cdot(\mu_\backslash,\mu_\tau) & \Omega \end{cases} \qquad (18)$$

where $\mathcal{H} \in \mathcal{L}(R^q, H^\backslash(\Omega))$ and $C:H^\backslash(\overline{\Omega}) \to R^q$ are a linear operator. If $n(b_\backslash - \alpha_\backslash)/(\beta_\backslash - \alpha_\backslash)$ and $m(b_\tau - \alpha_\tau)/(\beta_\tau - \alpha_\tau) \in \mathbb{N}$ for $n,m \in \{\backslash,\dots,J\}$, then the sensor $(b_i,\delta_{b_i})$ is not strategic for unstable sub-system of (16) and therefore the diffusion systems (16)-(17) are not detectable in $\Omega$ [20]. Then, the dynamical system (18) is not estimator for the diffusion systems (16)-(17) [1]. Here, we consider the region $\omega = ]\alpha_\backslash,\beta_\backslash[\times]\alpha_\tau,\beta_\tau[\subset]\cdot,a_\backslash[\times]\cdot,a_\tau[$ (Figure 3) and the dynamical system:

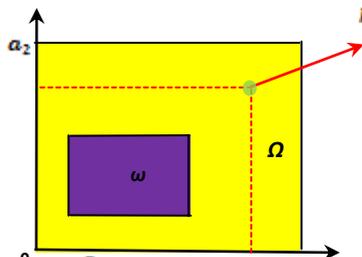

**Figure 3:** Domain $\Omega$, Region $\omega$, and position $b$ of the internal pointwise Sensor.

$$\begin{cases} \frac{\partial w}{\partial t}(\mu_\backslash,\mu_\tau,t) = \frac{\partial^\tau w}{\partial \mu_\backslash^\tau}(\mu_\backslash,\mu_\tau,t) + \frac{\partial^\tau w}{\partial \mu_\tau^\tau}(\mu_\backslash,\mu_\tau,t) & Q \\ \qquad\qquad + \mathcal{H}_\omega C\big(w(\mu_\backslash,\mu_\tau,t) - z(\mu_\backslash,\mu_\tau,t)\big) \\ \frac{\partial w}{\partial v}(\eta_\backslash,\eta_\tau,t) = \cdot & \Sigma \\ w(\mu_\backslash,\mu_\tau,\cdot) = w_\cdot(\mu_\backslash,\mu_\tau) & \Omega \end{cases} \qquad (19)$$

where $\mathcal{H}_\omega \in \mathcal{L}(R^q, H^\backslash(\omega))$. If $n(b_\backslash - \alpha_\backslash)/(\beta_\backslash - \alpha_\backslash)$ and $m(b_\tau - \alpha_\tau)/(\beta_\tau - \alpha_\tau) \in \mathbb{N}$ for every $n,m = \backslash,\dots,J$, then the sensor $(b_i,\delta_{b_i})$ is $\omega$-strategic for the unstable subsystem of (16) [18] and therefore the diffusion systems (16)-(17) are $\omega$-detectable, that means $\lim_{t\to\infty}\|z(\mu,t) - w(\mu,t)\|_{H^\backslash(\omega)} = \cdot$ [12]. Finally, the dynamical system (19) is $\omega$-estimator for the diffusion systems (16)-(17).

## 4. APPLICATION TO SENSORS LOCATIONS

In this section, we present an application of the above results to a two-dimensional diffusion system defined on $\Omega = ]\cdot, a_1[ \times ]\cdot, a_2[$ by the form

$$
\begin{cases}
\frac{\partial z}{\partial t}(\mu_1, \mu_2, t) = \frac{\partial^2 z}{\partial \mu_1^2}(\mu_1, \mu_2, t) + \frac{\partial^2 z}{\partial \mu_2^2}(\mu_1, \mu_2, t) & Q \\
\frac{\partial z}{\partial v}(\eta_1, \eta_2, t) = \cdot & \Sigma \\
z(\mu_1, \mu_2, t) = z_\cdot(\mu_1, \mu_2) & \Omega
\end{cases}
\tag{20}
$$

together with output function by (4) and (5) .Let $\omega = ]\alpha_1, \beta_1[ \times ]\alpha_2, \beta_2[$ be the considered region is subset of $]\cdot, a_1[ \times ]\cdot, a_2[$, $a_1 > \cdot, a_2 > \cdot$, and the eigenfunctions of system (20) are given by

$$
\varphi_{ij}(\mu_1, \mu_2) = \frac{1}{\sqrt{(\beta_1 - \alpha_1)(\beta_2 - \alpha_2)}} \cos i\pi \left(\frac{\mu_1 - \alpha_1}{\beta_1 - \alpha_1}\right) \cos j\pi \left(\frac{\mu_2 - \alpha_2}{\beta_2 - \alpha_2}\right)
\tag{21}
$$

associated with eigenvalues

$$
\lambda_{ij} = -\left(\frac{i^2}{(\beta_1 - \alpha_1)^2} + \frac{j^2}{(\beta_2 - \alpha_2)^2}\right) \pi^2.
\tag{22}
$$

The following results give information on the location of internal zone or pointwise $\omega$- strategic sensors.

### 4.1 Zone Sensor Case

Consider the systems (20)-(2) where the sensor supports $D$ is located in $\Omega$ .The output (2) can be written by the form

$$
y(t) = \int_D \frac{\partial x}{\partial v}(\mu_1, \mu_2, t) f(\mu_1, \mu_2) d\mu_1 d\mu_2
\tag{23}
$$

where $D \subset \Omega$ is the location of the internal zone sensor and $f \in H^1(D)$. In the case of (Figure 4), the eigenfunctions and the eigenvalues

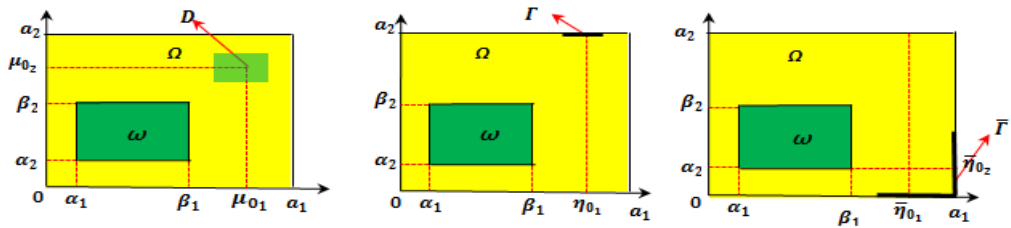

**Figure 4:** Domain $\Omega$, Sub-region $\omega$, and Locations of Internal (Boundary) Zone Sensor.

are given by (21)-(22). However, if we suppose that

$$
\frac{(\beta_1 - \alpha_1)^2}{(\beta_2 - \alpha_2)^2} \quad \text{is not even integers.}
\tag{24}
$$

then $r = 1$ and one sensor may be sufficient to achieve $\omega_{E_G}$-observability [12]. In this case, the dynamical system (13) is given by

$$
\begin{cases}
\frac{\partial w}{\partial t}(\mu_1, \mu_2, t) = \frac{\partial^2 w}{\partial \mu_1^2}(\mu_1, \mu_2, t) + \frac{\partial^2 w}{\partial \mu_2^2}(\mu_1, \mu_2, t) & Q \\
\qquad + \mathcal{H}_\omega(< z(\mu, t), f_i(\mu) > -Cw(\mu, t)) & \\
\frac{\partial w}{\partial v}(\eta_1, \eta_2, t) = \cdot & \Sigma \\
w(\mu_1, \mu_2, \cdot) = w_\cdot(\mu_1, \mu_2) & \Omega
\end{cases}
\tag{25}
$$

the measurement support is rectangular with

$$D = [\mu_\backslash, - I_\backslash, \mu_\backslash, + I_\intercal] \times [\mu_\intercal, - I_\backslash, \mu_\backslash, + I_\intercal] \in \Omega \tag{26}$$

If $f_\backslash$ is symmetric about $\mu_\backslash = \mu_{.\backslash}$, and $f_\intercal$ is symmetric about $\mu_\intercal = \mu_{.\intercal}$, then we have the following result.

**Proposition 4.1:** The systems (20)-(23) are not $\omega_{E_c}$-observable by the dynamical system (25) if

$\dfrac{I(\mu_{.\backslash} - \alpha_\backslash)}{(\beta_\backslash - \alpha_\backslash)}$ and $\dfrac{J(\mu_{.\intercal} - \alpha_\intercal)}{(\beta_\intercal - \alpha_\intercal)}$ , is integers for some $\ i = \backslash, \intercal, \dots, J$ .

In the case of boundary zone sensor, the system may be given by

$$\begin{cases} \dfrac{\partial w}{\partial t}(\mu_\backslash, \mu_\intercal, t) = \dfrac{\partial^\intercal w}{\partial \mu_\backslash^\intercal}(\mu_\backslash, \mu_\intercal, t) + \dfrac{\partial^\intercal w}{\partial \mu_\intercal^\intercal}(\mu_\backslash, \mu_\intercal, t) & Q \\ \qquad\qquad + \mathcal{H}_\omega (< z(\mu, t), f_i(\mu) >_{H^\backslash_{(\Gamma_i)}} - Cw(\mu, t)) & \\ \dfrac{\partial w}{\partial v}(\eta_\backslash, \eta_\intercal, t) = \cdot & \Sigma \\ w(\mu_\backslash, \mu_\intercal, \cdot) = w_{.}(\mu_\backslash, \mu_\intercal) & \Omega \end{cases} \tag{27}$$

and the output function (2) is given by

$$y(t) = \int_\Gamma \dfrac{\partial x}{\partial v}(\eta_\backslash, \eta_\intercal, t) f(\eta_\backslash, \eta_\intercal) d\eta_\backslash d\eta_\intercal \tag{28}$$

where $\Gamma \subset \partial\Omega$ is the support of the boundary sensor and $f \in H^\backslash(\Gamma)$. The sensor $(D, f)$ may be located on the boundary in $\Gamma = [\eta_{.\backslash}, - I, \eta_{.\backslash}, + I] \times \{a_\intercal\}$, then we have

**Proposition 4.2:** If $f_\backslash$ is symmetric about $\eta_\backslash = \eta_{.\backslash}$, then the systems (27)-(28) are not $\omega_{E_c}$-observable if $\dfrac{i(\eta_\backslash - \alpha_\backslash)}{(\beta_\backslash - \alpha_\backslash)}$ is integers for every $i, \backslash \le i \le J$.

When the sensor is located in $\bar{\Gamma} = [\bar\eta_{.\backslash}, - I_\backslash, a_\backslash] \times \{\cdot\} \cup \{a_\backslash\} \times [\cdot, \bar\eta_{.\intercal}, + I_\intercal] \subset \partial\Omega$ (see fig.3) we obtain the following result.

**Proposition 4.3:** suppose that the function $f|_{\Gamma_\backslash}$ is symmetric with respect to $\eta_\backslash = \bar\eta_{.\backslash}$, and the function $f|_{\Gamma_\backslash}$ is symmetric about $\eta_\intercal = \bar\eta_{.\intercal}$. Then the systems (27)-(28) are not $\omega_{E_c}$-observable if $\dfrac{(\eta_\backslash - \alpha_\backslash)}{(\beta_\backslash - \alpha_\backslash)}$ and $\dfrac{(\eta_\intercal - \alpha_\intercal)}{(\beta_\intercal - \alpha_\intercal)}$ is integers for every $i, \backslash \le i \le J$.

### 4.2 Pointwise Sensor Case

Let us consider the case of pointwise sensor located inside of $\Omega$. The system (20) is augmented with the following output function:

$$y(t) = \int_\Omega \dfrac{\partial z}{\partial v}(\mu_\backslash, \mu_\intercal, t) \delta(\mu_\backslash - b_\backslash, \mu_\intercal - b_\intercal) d\mu_\backslash d\mu_\intercal \tag{29}$$

where $b = (b_\backslash, b_\intercal)$ is the location of pointwise sensor as defined in (Figure 5)

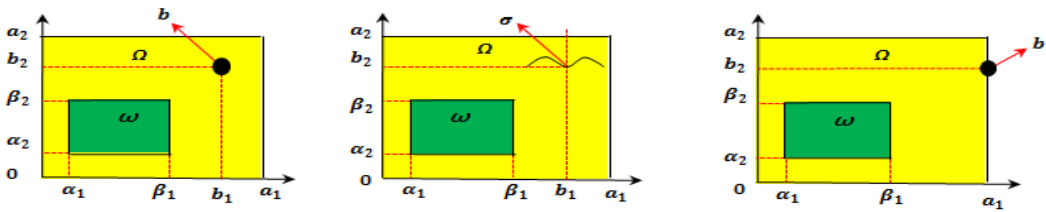

**Figure 5:** Rectangular Domain $\Omega$, and Locations $b$, $\sigma$ of Pointwise Sensor.

If $(\beta_1 - \alpha_1)^2/(\beta_2 - \alpha_2)^2 \notin Q$, then $r = 1$ and one sensor $(b, \delta_b)$ may be sufficient for $\omega_{E_G}$-observability. Then, the dynamical system is given by

$$\begin{cases} \dfrac{\partial w}{\partial t}(\mu_1,\mu_2,t) = \dfrac{\partial^2 w}{\partial \mu_1^2}(\mu_1,\mu_2,t) + \dfrac{\partial^2 w}{\partial \mu_2^2}(\mu_1,\mu_2,t) & Q \\ \qquad + \mathcal{H}_\omega\big(z(b_1,b_2,t) - Cw(\mu,t)\big) & \\ \dfrac{\partial w}{\partial v}(\eta_1,\eta_2,t) = 0 & \Sigma \\ w(\mu_1,\mu_2,\cdot) = w_{\cdot}(\mu_1,\mu_2) & \Omega \end{cases} \tag{30}$$

Thus, we obtain the following proposition:

**Proposition4.4:**The systems (20)-(29) are $\omega_{E_G}$-observable by the dynamical system (30) if $i(b_1 - \alpha_1)/(\beta_1 - \alpha_1)$ and $j(b_2 - \alpha_2)/(\beta_2 - \alpha_2)$, is not even integers for every $i, 1 \le i \le J$.

Consider  where the observation, in the case of the filament sensor, on the curve $\sigma = Im(\gamma)$ with $\gamma \in C^1(0,1)$ (Figure 4), then we have the following result

**Proposition 4.5:**If the observation recovered by filament sensor $(\sigma, \delta_\sigma)$ such that it is symmetric with respect to the line $b = (b_1, b_2)$, then the systems (20)-(29) are $\omega_{E_G}$-observable by (30) if $i(b_1 - \alpha_1)/(\beta_1 - \alpha_1)$ and $j(b_2 - \alpha_2)/(\beta_2 - \alpha_2)$, is not even integers for all $i = 1, \dots, J$.

And the case of the boundary pointwise sensor, the system (20) with Neumann boundary condition, so we can study the sensor $(b, \delta_b)$ is located on $\partial \Omega$ (figure 4). The output function is given by

$$y(t) = \int_\Omega \dfrac{\partial z}{\partial v}(\eta_1,\eta_2,t)\delta(\eta_1 - b_1, \eta_2 - b_2)d\eta_1 d\eta_2 \tag{31}$$

Then we can obtain.

**Proposition 4.6:** the systems (20)-(31) are $\omega_{E_G}$-observable by the dinamical system (30) if $i(b_1 - \alpha_1)/(\beta_1 - \alpha_1)$ and $j(b_2 - \alpha_2)/(\beta_2 - \alpha_2)$, is not even integers for every $i, 1 \le i \le J$.

**Remark 4.7:** From this work, we can extend to the following:

1. Case of mixed boundary conditions [1].

2. Case of disc domain $\Omega = (D, 1)$ and $\omega = (0, r_\omega)$ where $\omega \subset \Omega$

   and $0 < r_\omega < 1$ [10].

3. Case of one dimensional systems with pointwise or zone sensors as in

   ref.s [21, 5,8].

**المستخلص:** الهدف من هذه الورقة البحثية هو توسيع مفهوم قابلية المشاهدة العامة الاسية لمنطقة الى حالة مسالة شروط الحدودية لنيومان في نظام الانتشار . بشكل ادق ، في الانظمة الانتشار التوزيعية، نبين ان خصائص المجسات تسمح بإعادة بناء الحالة الاسية في منطقة جزئية من مجال الانظمة تحت الدراسة. اضافة على ذلك، يتم تقديم العديد من النتائج المثيرة للاهتمام والمرتبطة باختيار المجسات اعطيت ووضحت  في حالات محددة لكي يتم تحقيق قابلية المشاهدة العامة الاسية لمنطقة . واخيرا ، نبين ايضا ان هناك نظاما ديناميكيا لنظام الانتشار غير قابل للمشاهدة العامة الاسية بالمعنى العام، ولكن قد يكون غير قابل للمشاهدة الاسية العامة لمنطقة .

## 5. CONCLUSION

The regional exponential general observability in diffusion systems has been developed in this research and presented some sufficient condition of $\omega_{E_G}$-observability in connection with the Neumann conditions and strategic sensors. It allows us to avoid some "bad" sensor locations. Thus, the obtained results are applied to two-dimensional systems and various cases of sensors are considered and examined. Many questions are still opened, for example, the problem of finding the optimal sensor location ensuring such an objective. The result of regional exponential general observability concept of hyperbolic linear or semi-linear or nonlinear systems is under consideration.

**ACKNOWLEDGMENTS:** Our thanks in advance to the editors and experts for considering this paper to publish in this estimated journal and for their efforts

.